\documentclass{article}

\usepackage{arxiv}

\usepackage[utf8]{inputenc} % allow utf-8 input
\usepackage[T1]{fontenc}    % use 8-bit T1 fonts
\usepackage{hyperref}       % hyperlinks
\usepackage{url}            % simple URL typesetting
\usepackage{booktabs}       % professional-quality tables
\usepackage{amsfonts}       % blackboard math symbols
\usepackage{nicefrac}       % compact symbols for 1/2, etc.
\usepackage{microtype}      % microtypography
\usepackage{lipsum}		% Can be removed after putting your text content
\usepackage{graphicx}
\usepackage{subcaption}
\usepackage{natbib}
\usepackage{doi}
\usepackage{amssymb}
\usepackage{amsthm}
\usepackage{amsmath}
\newtheorem{remark}{Remark}
\newtheorem{definition}{Definition}
\newtheorem{theorem}{Theorem}

\title{Numerical methods that preserve a Lyapunov function for Ordinary Differential Equations}

\author{
	Yadira Hern{\'a}ndez-Solano \\
	Departmento de Matemática Aplicada\\
	Universidad de Málaga\\
	\texttt{yhdez@uma.es} \\
\And
\href{https://orcid.org/0000-0002-5158-5905}{\includegraphics[scale=0.06]{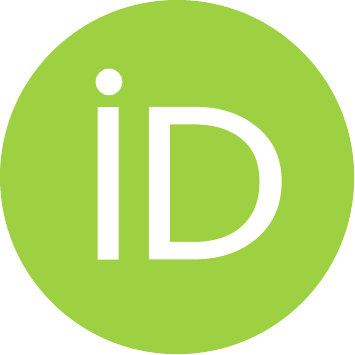}\hspace{1mm}Miguel Atencia} \\
	Departmento de Matemática Aplicada\\
	Universidad de Málaga\\
	\texttt{matencia@uma.es} \\
}

\begin{document}
\maketitle

\begin{abstract}
The paper studies numerical methods that preserve a Lyapunov function of a dynamical system, i.e. numerical approximations whose energy decreases, just like in the original differential equation. With this aim, a discrete gradient method is implemented for numerical integration of a system of ordinary differential equations. In principle, this procedure yields first order methods, but the analysis paves the way to the design of higher-order methods. As a case in point, the proposed method is applied to the Duffing equation without external forcing, considering that in this case, preserving the Lyapunov function is more important than accuracy of particular trajectories. Results are validated by means of numerical experiments, where the discrete gradient method is compared to standard Runge-Kutta methods. As predicted by the theory, discrete gradient methods preserve the Lyapunov function, whereas conventional methods fail to do so, since either periodic solutions appear or the energy does not decrease. Besides, the discrete gradient method outperforms conventional schemes when these do preserve the Lyapunov function, in terms of computational cost, thus the proposed method is promising.
\end{abstract}

% keywords can be removed
\keywords{Geometric Numerical Integration \and Dynamical Systems \and Lyapunov Function \and Stability \and Numerical Methods \and Ordinary Differential Equations \and Discrete Gradient}

\section{Introduction}
The main aim of this paper is the study of numerical methods that preserve a Lyapunov function of a gradient dynamical system. The solutions, or integral curves, of a gradient system follow trajectories that are tangent to a scalar function of the states, which is usually known as the Lyapunov function of the system. The flow of a gradient system has a rather simple qualitative behaviour, e.g. all isolated minima of the Lyapunov function are asymptotically stable equilibria of the system \cite{hirsch}. The Lyapunov function has the remarkable property that it is decreasing along trajectories of the dynamical system. Gradient systems are pervasive, both as models of physical systems and as representations of mathematical algorithms. For example, an ideal pendulum is a conservative system, namely the energy is a constant magnitude, but every actual mechanical system dissipates energy due to friction, until all potential and kinetic energies vanish, thus energy acts as Lyapunov function of the system. Remarkably, many mathematical algorithms are formulated in continuous time whose operation is based on the existence of the Lyapunov function, for example, in the fields of optimization, estimation, and control \cite{absil_continuous_2004,slotine}.

Numerical methods for the integration of Ordinary Differential Equations (ODEs) constitute a well established field \cite{hairer}, and methods that provide rather accurate solutions for a wide variety of problems have long been known. However, no matter how small the approximation error of a numerical method is, it can lead to a solution that does not portray the qualitative features of the continuous model, when the integration extends through long time periods. A classic example is the Kepler problem \cite{arnold_mathematical_1997}, whose approximate solutions by conventional numerical methods do not respect the elliptical orbits describing the motion of the planets, as established by Kepler's first law. The inability of basic numerical methods to reflect crucial qualitative properties of dynamical systems, led to the development of a new approach, namely \emph{Geometric Numerical Integration} \cite{hairer2}, which is an active line of research that links the methodology of dynamical systems analysis to the design of numerical methods~\cite{stu-hum} that preserve the qualitative properties of the continuous system. In this regard, the main objective is to consider the qualitative characteristics of the trajectories of the dynamical system, e.g. energy decreasing, stability, conservation of the Hamiltonian, among others. The task would then be the design of numerical methods, so that the discrete trajectories of the method  have the same properties as the exact solutions.

Within the field of Geometric Numerical Integration, there exists a substantial number of results concerning the study of Hamiltonian systems \cite{sanz-serna_symplectic_2008}, among which symplectic and projection methods can be mentioned. However, when it comes to the conservation of the Lyapunov function of a gradient dynamical system,  the choice is limited---to the best of our knowledge---to three categories: discrete gradient methods \cite{quispel}, projection methods~\cite{calvo_projection_2010}, and particular instances of Runge-Kutta methods \cite{hairer_energy-diminishing_2013}. The inattention to stability issues is striking, since dynamic analysis of ODEs is far from new in the field of numerical analysis. Indeed, the concept of A-stability \cite{iserles_first_2009} amounts to the preservation of the stability of the solution of linear scalar equations as test systems. In this regard the conservation of the Lyapunov function can be viewed as a generalization of the concept of A-stability in a nonlinear context.

Discrete gradient methods \cite{quispel} yield integrators for ODEs, based upon the fact that the equation of a gradient system can be written in linear-gradient form, i.e. as the product of a definite-negative matrix by the gradient of a Lyapunov function. Then, discrete gradient methods can be stated with a simple rationale: define an approximation of the definite-negative matrix and a discrete gradient, which has similar properties to that of the gradient of the Lyapunov function. By construction, these methods lead to an implicitly defined map that, when considered as a discrete dynamical system, preserves that Lyapunov function of the continuous system. The development of discrete gradient methods is limited and examples of systematic application to real systems are hardly found in the literature, as far as we know. In previous work, we explored the application of discrete gradient methods to a particular system, namely Hopfield neural networks \cite{iwann_2013,Hernandez-Solano2015}, which are computational methods used for optimization. For its part, projection methods~\cite{calvo_projection_2010} inherit the design of analogous methods for Hamiltonian systems, which are based upon projecting the approximate solution onto the manifold that the trajectories of the exact solutions lie in. Although the formulation of these methods is explicit in principle, they require to solve the nonlinear equation that defines the projection at each time step. Finally, the application of Runge-Kutta methods to gradient systems~\cite{hairer_energy-diminishing_2013} led to proving that some Radau implicit methods, originally proposed for stiff and Hamiltonian systems, are also able to preserve the Lyapunov function under certain conditions and restrictions on the step size. Both projection methods and implicit Runge-Kutta integrators rely on non-constructive theorems, so they cannot guarantee the preservation of the Lyapunov function unless some ad-hoc adjustment of the step size is performed. In summary, although these methods are promising, their implementation is complicated and can lead to a substantial computational cost, so that they are not suitable for all situations. It must be emphasized that, rather than advocating against other techniques, our results encourage further attention to discrete gradient methods, at least for particular applications. Nonetheless, some considerations on future lines of research for a comparative assessments of all these methods are made in the conclusions.

After a review of the background about discrete gradient methods in Section \ref{sec:methods}, the contribution of this paper begins in Section \ref{sec:dg}, where we describe the methodology of implementation of discrete gradient methods, analysing the order of the obtained method and illustrating its main properties by means of simple examples. Then, in Section \ref{sec:experiment} we present some systematic numerical experiments showing the performance of the proposed technique, and comparing its performance to standard Runge-Kutta methods. As a result, some favourable properties of the obtained method are brought to light. Finally, some conclusions and lines for further research are stated in Section~\ref{sec:conclusion}.

\section{Numerical methods that preserve the Lyapunov function}
\label{sec:methods}
In this section, we define and discuss key aspects of discrete gradient methods, after establishing the definitions that will be used along the paper.

\subsection{Gradient systems}
First of all, we establish the notation for the dynamical system that must be dealt with, which is a finite-dimensional initial value problem (IVP), i.e. a system of ODEs with initial values:
\begin{equation}
\label{eq1}
    \frac{d\vec{y}}{dt}=\vec{f}(\vec{y})=
\begin{pmatrix}
f_{1}(\vec{y}) \\
\vdots \\
f_{n}(\vec{y})
\end{pmatrix},
\qquad \vec{y}(t_{0})=\vec{y}_{0}\in \mathbb{R}^{n}
\end{equation}
Since we do not pursue existence and uniqueness issues, we take for granted all needed smoothness assumptions. The systems of interest are those that possess---at least---one asymptotically stable equilibrium (see e.g. \cite{hirsch,khalil} for definitions of stability concepts). An equilibrium or fixed point $\vec{y}^*$ fulfills $\vec{f}(\vec{y}^*)=\vec{0}$, thus a trajectory that starts at $\vec{y}_{0}=\vec{y}^*$ is the trivial trajectory $\vec{y}(t)=\vec{y}^*$. The statement that $\vec{y}^*$ is asymptotically stable amounts to saying that all trajectories $\vec{y}(t)$ that start in a certain neighbourhood $\mathcal{B}$ with $\vec{y}^* \in \mathcal{B}$ converge towards $\vec{y}^*$, i.e.: $\displaystyle\lim_{t\to\infty} \vec{y}(t)=\vec{y}^*$ if $\vec{y}_{0}\in\mathcal{B}$.

One of the aims of qualitative analysis of ODEs is proving that an equilibrium is stable without computing the explicit solution, which can be accomplished by finding a suitable Lyapunov function $V$:
\begin{definition}
\label{def:lyapunov}
Given the system in Equation \eqref{eq1}, the function $V\in C^{1}(\mathbb{R}^{n}, \mathbb{R})$ is a Lyapunov function for the equilibrium $\vec{y}^*$ if the following conditions hold in a neighbourhood~$\mathcal{B}$ such that $\vec{y}^* \in \mathcal{B}$:\\
\begin{itemize}
  \item [a)] $V(\vec{y}^*)=0$ and $V(\vec{y})>0$ if $\vec{y}\neq\vec{y}^*$. \\
  \item [b)] $\displaystyle\frac{d}{dt}{V(\vec{y}(t))} < 0$  for all $\vec{y} \in
                 \mathcal{B}-\left\{\vec{y}^*\right\}$. \\
\end{itemize}
\end{definition}
By the chain rule, this definition is equivalent to stating the following relation of the gradient and the ODE:
$\displaystyle\frac{dV}{dt}=\nabla V(\vec{y}) \cdot \vec{f}(\vec{y})\le 0$ with $V$ bounded below. Note that it is obvious that $\displaystyle\frac{dV}{dt}=0$ at an equilibrium $\vec{y}^*$ since $ \vec{f}(\vec{y}^*)=\vec{0}$. The existence of a Lyapunov function characterizes the stability of an equilibrium \cite{khalil}:
\begin{theorem}
Let $\vec{y}^*$ be an equilibrium point of the system in Equation \eqref{eq1} and $V$ a Lyapunov function in a neighborhood $\mathcal{B}$ of $\vec{y}^*$. Then, $\vec{y}^*$ is asymptotically stable.
\end{theorem}

A Lyapunov function is often called the \emph{energy} of the system, by analogy with dissipative physical systems where the energy decreases, thus it may be used as a Lyapunov function. Rigorously speaking, the definition of Lyapunov function does not require the inequality in condition b) of Definition~\ref{def:lyapunov} to be strict and, when the inequality is strict, we should specify that the system has \emph{strict} Lyapunov function. In this paper we always assume that the Lyapunov function is strict, so we do not make this distinction. Likewise we loosely refer to a stable point, dropping the assumed precision that such stability is asymptotic. It is worth remarking that, although converse theorems guarantee the existence of a Lyapunov function when a stable equilibrium exists, there is no general method for finding the explicit expression of a Lyapunov function. In this paper, we assume that a Lyapunov function is explicitly known.

The main aim of this paper is to find numerical methods that preserve the Lyapunov function of a system given by Equation \eqref{eq1}. Formally, we construct a discrete dynamical system defined by a time-stepping formula $\vec{z}=\varphi_h\left(\vec{y}\right)$ such that $\vec{z}$ is a suitable approximation of $\vec{y}(t+h)$ if $\vec{y}$ is an approximation of $\vec{y}(t)$. The required preservation of the Lyapunov function $V$ is subsumed by the condition $V(\vec{z})<V(\vec{y})$ as long as $\vec{z}\neq\vec{y}$, which is the discrete counterpart of condition b) in Definition \ref{def:lyapunov}: both inequalities express that the Lyapunov function decreases through time, either in a discrete or a continuous setting. It will also be of interest to determine if the time-stepping scheme produces a sequence that converges to some stable equilibrium of the original system, thus reproducing asymptotic stability.

\subsection{Discrete gradient methods}
The history of stability preserving methods can be traced back at least three decades, to the seminal paper \cite{Hairer1990} and, later, the book \cite{stu-hum}. It is thus well known that numerical methods may destroy the structural properties of the original ODE, but note that there is a hierarchy of how subtle can this effect be. On the one hand, an equilibrium may cease to be fixed point of the discrete method, or it may become an unstable equilibrium. These spurious solutions can easily be detected by a (more or less) straightforward analysis of the method, including linearization around the equilibrium. More importantly, there are established criteria to construct (local) stability-preserving numerical methods. A much more severe problem arises when the equilibrium is still locally asymptotically stable, but the numerical method fails to decrease the Lyapunov function or, in other words, the basins of attraction change. This alteration of geometrical properties has a global nature, hence its study is notoriously difficult. Discrete gradient methods guarantee that the Lyapunov function of the ODE decreases along sequences of points obtained by the numerical method so that, at least from the point of view of energy minimization, the geometric structure is preserved.

The rationale behind discrete gradient methods is a rather simple idea, namely to replace the derivative of the Lyapunov function by a finite increment. This idea is useful for discretizing the system, because the ODE and the Lyapunov function are related: every ODE as in Equation \eqref{eq1} for which a Lyapunov function $V$ is known, can be rewritten in \emph{linear-gradient} form~\cite{quispel}:
\begin{equation}
\label{GL}
    \frac{d\vec{y}}{dt}=L(\vec{y})\,\nabla V(\vec{y})
\end{equation}
where $L$ is a negative-definite matrix and both $L$ and $V$ are continuously differentiable.
Incidentally, it is worth mentioning that this decomposition is not unique, and the different ways to write $L(\vec{y})$ can be regarded as different metric structures \cite{barta_every_2011}.

\begin{remark}
Care must be taken when negative-definiteness is considered for non-symmetric matrices, since in this case negative eigenvalues of $L$ do not guarantee the intended relation $v^\top L v<0$ for any vector $v\neq 0$. Let us thus emphasize that, along the paper, a matrix $L$ is negative-definite if its \emph{symmetric part} $L+L^\top$ is.
\end{remark}

After rewriting the ODE in Equation \eqref{eq1} in linear-gradient form, a discrete gradient method results from the choice of discrete analogs to the matrix $L$ and the gradient $\nabla V$:
\begin{definition}
Given a differentiable function $V\in C^{1}(\mathbb{R}^{n}, \mathbb{R})$, the function $\overline{\nabla}V \in C^{1}(\mathbb{R}^{2\,n}, \mathbb{R}^n)$ is a \emph{discrete gradient} of $V$ if it satisfies:
\begin{equation}
\label{gd-conditions}
\begin{aligned}
\overline{\nabla }V\left(\vec{y},\vec{z}\right) \cdot \left(\vec{z}-\vec{y}\right)
&=V\left(\vec{z}\right)-V\left(\vec{y}\right)
\\
\overline{\nabla }V\left(\vec{y},\vec{y}\right)&=\nabla V\left(\vec{y}\right)
\end{aligned}
\end{equation}
\end{definition}
In fact, the second condition is implied by the first in the differentiable case \cite{Eidnes2022}, but we include it anyway to emphasize consistency.

\begin{definition}
A \emph{discrete gradient method} is a time-advancing numerical scheme defined by
\begin{equation}
\label{MI}
    \frac{\vec{z}-\vec{y}}{h}=\widetilde{L}(\vec{y},\vec{z},h)\;\overline{\nabla}V(\vec{y},\vec{z})
\end{equation}
where $\overline{\nabla}V$ is a discrete gradient of $V$ and the matrix $\widetilde{L}(\vec{y},\vec{z},h)$ of continuously differentiable functions is negative definite and satisfies the consistency condition
\begin{equation}
\label{l-condition}
\widetilde{L}\left(\vec{y},\vec{y},0\right)=L\left(\vec{y}\right)
\end{equation}
\end{definition}

The aim of a discrete gradient method is to compute $\vec{z}\approx\vec{y}(t+h)$ from the previous step $\vec{y}\approx\vec{y}(t)$ so the sequence $\vec{y}(t)$ is an approximation of the solution of the system given by Equation \eqref{GL}. It is trivial to prove that a discrete gradient method is consistent, as a consequence of the requirements on  $\widetilde{L}$ and  $\overline{\nabla}V$. Remarkably, the methods given by Equation \eqref{MI} are implicit, at least in principle, since the next step $\vec{z}$ appears in the right-hand side of the formula.

\section{Construction and analysis of Discrete Gradient Methods}
\label{sec:dg}
Once the parameters $\widetilde{L}$ and $\overline{\nabla }V$ have been set, a particular instance of discrete gradient method results by substituting this parameter choice into Equation \eqref{MI}. This is a critical design process, since there is a wide range of choices that are compatible with the conditions given by Equations \eqref{gd-conditions} and \eqref{l-condition}. Regarding the matrix $\widetilde{L}$, the trivial choice  $\widetilde{L}(\vec{y},\vec{z},h) = L\left(\vec{y}\right)$ is possible, where the dependence on the next step $\vec{z}$ is neglected. A less radical simplification results when dismissing the step size $h$ in the definition of  $\widetilde{L}$. We adopt this latter assumption throughout this paper, so we often write $\widetilde{L}(\vec{y},\vec{z})$ for this matrix. With regard to the discrete gradient, there is a single discrete gradient for one-dimensional systems,  and it is given by:
\begin{equation}
\overline{\nabla}V\left(\vec{y},\vec{z}\right)=\frac{V(\vec{z})-V(\vec{y})}{\vec{z}-\vec{y}}
\end{equation}
However in higher dimensions a wide variety of discrete gradients exist (see~\cite{quispel} and references therein for several examples). In this paper we will focus on the \emph{coordinate increment discrete gradient}, also called Itoh-Abe discrete gradient~\cite{ITOH198885}, since it is easier to implement computationally. We assume an ordering $y_{1}, y_{2},\ldots, y_{n}$ of the coordinates of the vector $\vec{y}\in\mathbb{R}^n$ and define the function $\overline{\nabla}V$ as:
\begin{equation}
\label{gd}
\overline{\nabla}V\left(\vec{y},\vec{z}\right)=
\begin{pmatrix}
\displaystyle\frac{V\left(z_{1},y_{2},...,y_{n}\right)-V\left(
y_{1},y_{2},...,y_{n}\right)
}{z_{1}-y_{1}}\\ \\
\displaystyle\frac{V\left(z_{1},z_{2},y_{3},...,y_{n}\right)-V\left(
z_{1},y_{2},...,y_{n}\right)
}{z_{2}-y_{2}}\\ \vdots \\
\displaystyle\frac{V\left(z_{1},...,z_{n-2},z_{n-1},y_{n}\right)-V\left(
z_{1},...,z_{n-2},y_{n-1},y_{n}\right)
}{z_{n-1}-y_{n-1}}\\\\
\displaystyle\frac{V\left(z_{1},...,z_{n}\right)-V\left(
z_{1},...,z_{n-1},y_{n}\right) }{z_{n}-y_{n}}
\end{pmatrix}
\end{equation}
The coordinate increment discrete gradient can be interpreted as a piecewise linear path joining $\vec{y}$ and $\vec{z}$, each piece parallel to one of the coordinate axes, rather than along the straight segment $\vec{y}-\vec{z}$.

In the rest of this section, we undertake a study of discrete gradient methods, first by a preliminary order analysis, then by constructing different methods for simple scalar systems (this methodology is inspired by \cite{ramos_piecewise-linearized_1997}) and observing that a suitable choice of the matrix  $\widetilde{L}$ allows in some cases for rewriting the method in explicit form.

\subsection{Order analysis}
The order of the obtained numerical method can be studied by the usual systematic procedure~\cite{hairer}: comparing the Taylor series expansion around $h=0$ of both the exact solution of the system of differential equations and the approximate solution obtained by the numerical method. Note that the discrete gradient method is consistent by construction \cite{quispel} so it achieves at least order one, i.e. the error after a single step is given by $\vec{y}(t+h)-\vec{z} = C \, h^{2}+O(h^3)$, where $C$ is the error constant of the method. A straightforward---but tedious---computation yields  the error constant of the second order term:
\begin{equation}
\label{consterror_gd}
C_{GD}=\left( \dfrac{1}{2}\, J^{c} -\left.J^{d}\right|_{h=0} \right) \, f(\vec{y})
\end{equation}
where $J\,^{c}(\vec{y})$ is the Jacobian matrix of $f$ at $\vec{y}$:
\begin{equation}
\label{Jcont}
    J\,^{c}=\frac{\partial f}{\partial \vec{y}}=\left[\frac{\partial f_{i}}{\partial y_{j}}\right]_{ij} =
\left[\frac{\partial \left(L(\vec{y}) \cdot \nabla V(\vec{y})\right)_i}{\partial y_{j}}\right]_{ij}
\qquad    i,j=1,\ldots,n
\end{equation}
and  $J^d$ is the Jacobian of $\widetilde{L}\cdot \overline{\nabla}V$, i.e.:
\begin{equation}
\label{jacobianodiscreto}
    J\,^{d}= \frac{\partial\left(\widetilde{L}\cdot \overline{\nabla}V\right)}{\partial
    \textbf{z}}=\left[\frac{\partial\left(\widetilde{L}\cdot \overline{\nabla}V\right)_{i}}{\partial
    \textbf{z}_{j}}\right]_{ij} \quad i,j=1,\ldots,n
\end{equation}
so that the condition $J^{c} =2\,\left.J^{d}\right|_{h=0}$ would ensure that the obtained discrete gradient method is second order. In principle, a suitable choice of parameters $\widetilde{L}$ and $\overline{\nabla}V$ could lead to a higher-order method. When this paper was already in preparation, a systematic analysis of discrete gradient methods has been published \cite{Eidnes2022}, although in the somewhat different context of Hamiltonian systems. Adapting this framework to gradient-like systems is an interesting avenue for future research.  Nevertheless it must be emphasized that the search for higher accuracy without any other consideration, defeats the purpose of structure preserving methods. In this paper we will not further pursue the analysis of order and error, focusing on the preservation of the Lyapunov function and stability.

\subsection{The scalar linear ODE}
For the purpose of illustration, in this section we show the mechanism of obtaining a discrete gradient method as described above. As a case in point, consider the scalar linear homogeneous  ODE:
\begin{equation}
\label{P1}
\frac{d y}{dt}=-a\,y\, , \qquad y(0)=y_0
\end{equation}
with $a>0$. By direct integration, it is straightforward to compute the analytical solution $y(t)=y_0 \, e^{-a\,t}$, which shows that the origin is asymptotically stable whenever $a>0$, since $\displaystyle\lim_{t\to\infty}y(t)=0$. We can also state that $V = \dfrac{1}{2}\,y^{2}$ is a Lyapunov function for this system because:
\begin{equation}
\alpha(y)=\dfrac{dV}{dt}=\dfrac{dV}{dy}\dfrac{dy}{dt}=y\,(-a\,y)=-a\,y^2<0
\end{equation}
for all $y\neq 0$. In order to construct a discrete gradient method, the equation is cast into linear-gradient form, thus obtaining the definitions $L(y)=-a$, $\nabla V=y$. Therefore the discrete gradient is:
\begin{equation}
\label{GDP1dim1}
    \overline{\nabla}V(y,z)=\frac{V(z)-V(y)}{z-y}=\frac{1}{2}\,\frac{z^2-y^2}{z-y} =\frac{1}{2}\,(z+y)
\end{equation}
and, with the trivial choice $\widetilde{L}=L=-a$, the discrete gradient method results:
\begin{equation}
  z=y+h\,\widetilde{L}(y,z) \,  \overline{\nabla}V(y,z)=y+\frac{h}{2}\,(-a\,z-a\,y)=y+\frac{h}{2}\,\left(f(z)+f(y)\right)
\end{equation}
Now it is obvious that in this case the discrete gradient method turns out to be simply the trapezoidal rule, which is a second-order method. The fact that the trapezoidal rule preserves the stability of scalar linear ODEs for \emph{any} step size $h$ is already explained by the classical theory of numerical methods for stiff systems, since it is well-known that the trapezoidal rule is A-stable, thus nothing new seems to be provided by the proposal of discrete gradient methods. However, the point is that the choice of the matrix  $\widetilde{L}$ is not unique, so a different definition  $\widetilde{L}$, possibly depending on $\vec{z}$ and $h$, would lead to a different method. In addition, if we are not interested in preserving a particular Lyapunov function,  but only the qualitative stability of the system, we could choose a different Lyapunov function, thus leading to a different discrete gradient method. 

\subsection{The logistic equation}
Consider next the IVP given by the generalization of the usually called logistic differential equation:
\begin{equation}
\label{P4}
\frac{dy}{dt}=a \, y\, (1-y)\, , \qquad y(0)=y_0
\end{equation}
By straightforward integration, the exact solution can be computed:
\begin{equation}
\label{exactP4}
y(t)=\frac{1}{1+\left(\dfrac{1}{y_0}-1\right)\, e^{-at}}
\end{equation}
for any initial condition $y_0\ne 0$, whereas the trivial solution $y(t)=0$ involves a fixed point. We also choose $y_0>0$ to avoid the need to consider unbounded solutions. There are several ways to check that the equilibrium $y^*=1$ is asymptotically stable, e.g. the Jacobian of the ODE given by Equation \eqref{P4} is negative at $y=1$ or the limit when $t\to\infty$ of the exact solution given by Equation \eqref{exactP4} is $1$. 

The construction of  a discrete gradient method as in Equation \eqref{MI} requires, first, writing the system in linear-gradient form from the knowledge of a Lyapunov function $V$; and then choosing the method parameters, $\widetilde{L}(y,z,h)$ and $\overline{\nabla}V(y,z)$, while fulfilling the conditions that guarantee the consistency of the method. Interestingly, even such a simple example as the logistic ODE can lead to completely different  discrete gradient methods.

Firstly, observe that the function $V=\displaystyle\frac{1}{2}(1-y)^2$ fulfils the conditions required by Definition \eqref{def:lyapunov} to be a Lyapunov function. In particular, its time derivative is:
\begin{equation}
\label{eq:V-decreasing}
\frac{dV}{dt}=\nabla V\cdot f=-(1-y)\,  ay\, (1-y)=-ay\,(1-y)^2<0
\end{equation}
whenever $y>0,y \neq 1$. Therefore $V$ is a Lyapunov function of Equation \eqref{P4} at $y^*=1$ that is valid for any initial value $y_0>0$. Then, the ODE can be rewritten in linear-gradient form as in Equation \eqref{GL} by defining $L(y)=-a\,y$, so that the system is expressed as:
\begin{equation}
\label{GLP4}
    \frac{dy}{dt}=a\,y\;(1-y)=L(y)\,\nabla V
\end{equation}
with $L(y)$ negative-definite for $y>0$, as required. Then, the discrete gradient is defined by the unique choice existing in the scalar case:
\begin{equation}
\label{GDlogistic1}
\begin{aligned}
    \overline{\nabla}V(y,z)&=\frac{V(z)-V(y)}{z-y}=\frac{1}{2}\,\frac{(1-z)^2-(1-y)^2}{z-y} =\\
&=\frac{1}{2}\,\frac{-2\,(z-y)+(z^2-y^2)}{z-y}=\frac{-2+(z+y)}{2} \\
&= - \left(1-  \frac{z+y}{2}  \right)
\end{aligned}
\end{equation}
The last equality of Equation \eqref{GDlogistic1} has been included to point out a plausible interpretation of the discrete gradient as a sort of \emph{midpoint} gradient, since it is identical to the gradient of $V$, replacing the variable $y$ with the average $\dfrac{z+y}{2}$. With regard to the choice of $\widetilde{L}(y,z,h)$, there are several consistent options. For simplicity, we adopt the trivial setting $\widetilde{L}=L$. Therefore, if we substitute the chosen parameters  in Equation~\eqref{MI}, the method is obtained:
\begin{equation}
\begin{aligned}
    z&=y+h\,\widetilde{L}(y,z,h)\cdot\overline{\nabla}V(y,z)=y - a\,h \, y \, \left(\frac{-2+z+y}{2}\right) \\
    &= y  + a\, h \, y -  \frac{a\,h}{2} \,y\,z - \frac{a\,h}{2} \,y^2
\end{aligned}
\end{equation}
which after straightforward algebra yields an explicit expression for $z$:
\begin{equation}
\label{logistic-explicit}
    z=\frac{\left(1+a\,h-\dfrac{a\,h}{2}\,y\right)\;y}{1+\dfrac{a\,h}{2} \; y}
\end{equation}
In this particular case, the choice of $\widetilde{L}$ has allowed for obtaining an explicit method. However, the procedure has some generality, at least restricted to one-dimensional ODEs: it can be proved that if the Lyapunov function $V$ is quadratic and the matrix $\widetilde{L}$ is trivially set to $\widetilde{L}=L$, the discrete gradient method can be cast into explicit form. 

\begin{remark}[Relation to known methods]
Note that apparently Equation \eqref{logistic-explicit} cannot be derived as a conventional Runge-Kutta method (although proving this in general would require some work). In contrast, the \emph{nonlocal} substitution $y^2 \rightarrow yz$ and the use of the discrete gradient remind of nonstandard finite difference schemes \cite{Mickens2005}, while providing a systematic methodology for their construction. 
\end{remark}

%%%%%%%%%%%%%%%%%%%%%%%%%%%%%%%%%%%%%%%%%%%%%%%%%%%%%%%
Consider now the function $V=\displaystyle-\frac{1}{2}y^2+\frac{1}{3}y^3$ as a candidate for Lyapunov function of the same system, and observe that it fulfils the conditions required by Definition~\eqref{def:lyapunov}. In particular, the time derivative is:
\begin{equation}
\frac{dV}{dt}=\nabla V\cdot f=-y\,(1-y)\,  ay\, (1-y)=-ay^2\,(1-y)^2<0
\end{equation}
whenever $y\neq 0,1$. Therefore $V$ is a Lyapunov function of Equation \eqref{P4} for the stable  equilibrium point $y^*=1$, that is valid for any initial value $y_0> 0$. Then, the ODE can be rewritten in linear-gradient form as in Equation \eqref{GL} by defining $L(y)=-a$ so that the linear gradient form $
    \frac{dy}{dt}=ay\;(1-y)=L(y)\,\nabla V$
holds too with these new parameters and $L(y)$ is negative-definite, as required. The one-dimensional discrete gradient has the same form as before, but the Lyapunov function $V$ is different to begin with, leading to:
\begin{equation}
\label{GDP4dim1}
\begin{aligned}
    \overline{\nabla}V(y,z)&=\frac{V(z)-V(y)}{z-y}=\dfrac{-\dfrac{1}{2}(z^2-y^2)+\dfrac{1}{3}(z^3-y^3)}{z-y}=\\
&=-\frac{1}{2}(z+y)+\frac{1}{3}(z^2+zy+y^2)
\end{aligned}
\end{equation}
With regard to the choice of $\widetilde{L}(y,z,h)$, for simplicity we again adopt the trivial setting $\widetilde{L}=L$. Therefore, if we substitute the chosen parameters  in Equation~\eqref{MI}, the new method is obtained:
\begin{equation}
    z=y+h\,\widetilde{L}(y,z,h)\cdot\overline{\nabla}V(y,z)=y - ah \,\left[-\dfrac{1}{2}(z+y)+\dfrac{1}{3}(z^2+zy+y^2)\right]
\end{equation}
In this case we obtain an implicit method. In order to apply Newton's method to obtain the solutions, we can rewrite the method as a function of $z$ as shown below:
\begin{equation}
\label{logistic-implicit}
    F(z)=\left(\dfrac{1}{3}ah\right)z^2+\left(1-\dfrac{1}{2}ah+\dfrac{1}{3}ahy\right)z+\left(-1-\dfrac{1}{2}ah+\dfrac{1}{3}ahy\right)\,y=0
\end{equation}

\begin{figure}
\begin{subfigure}[b]{0.5\linewidth}
    \centering
    \includegraphics[width=\textwidth]{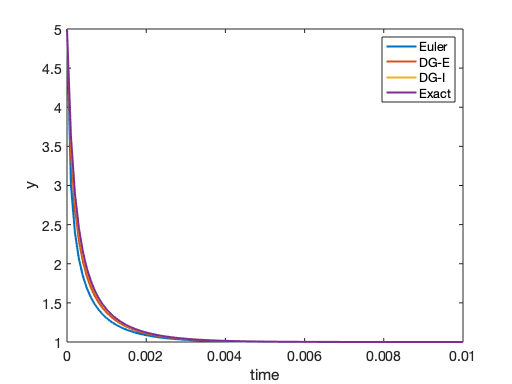}
    \caption{Trajectories for $h=10^{-4}$.}
    \label{subfig:logistica_todos}
\end{subfigure}
\hfill
\begin{subfigure}[b]{0.5\linewidth}
    \centering
    \includegraphics[width=\textwidth]{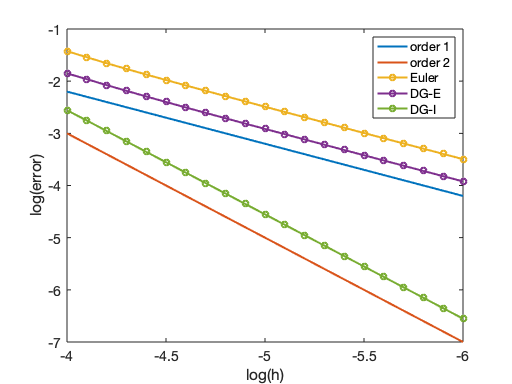}
    \caption{Error estimation.}
    \label{subfig:logistica_errors}
\end{subfigure}
\\\\
\begin{subfigure}[b]{0.5\linewidth}
    \centering
    \includegraphics[width=\textwidth]{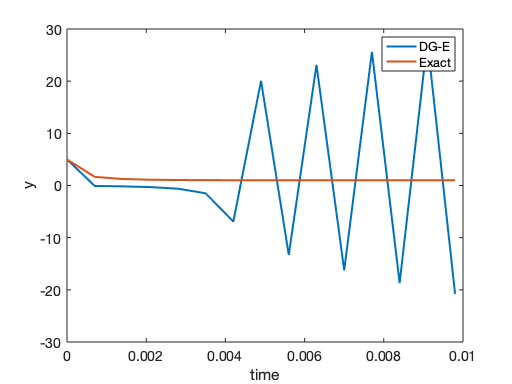}
    \caption{Trajectory for DG-E with $h=7\cdot 10^{-4}$.}
    \label{subfig:logistica_expl}
\end{subfigure}
\hfill
\begin{subfigure}[b]{0.5\linewidth}
    \centering
    \includegraphics[width=\textwidth]{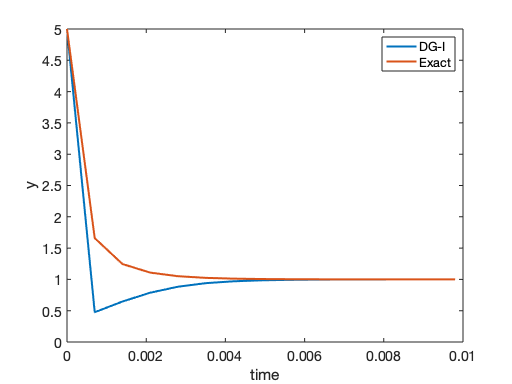}
    \caption{Trajectory for DG-I with $h=7\cdot 10^{-4}$.}
    \label{subfig:logistica_impl}
\end{subfigure}
    \caption{Solutions for the logistic equation obtained by the Euler method, the explicit method in Equation \eqref{logistic-explicit} (DG-E), the implicit method in Equation \eqref{logistic-implicit} (DG-I), and the exact solution.}
    \label{fig:logistica}
\end{figure}

We have implemented the explicit discrete gradient method (DG-E) given by Equation \eqref{logistic-explicit}, and the implicit scheme (DG-I) from Equation \eqref{logistic-implicit}. Both are applied to the same logistic ODE, choosing the parameter as $a=1000$ and the initial value $y_0=5$. The resulting trajectories are shown in Figure \ref{fig:logistica} for different values of the step size $h$. When $h$ is small enough all methods provide qualitatively correct solutions, as shown in Figure \ref{subfig:logistica_todos}. Besides both discrete gradient methods derived above, the Euler rule has been included for comparison. In order to have a glimpse at the approximation accuracy achieved by each method, the global error has been computed by subtracting the discrete sequence from the exact solution and averaging over all the computed steps. The obtained results for 20 different values of the step size in the interval $h \in [10^{-6},10^{-4}]$ are shown in Figure \ref{subfig:logistica_errors}, in logarithmic scale. Two straight lines with slopes 1 and 2 are added to ease the comparison. It is clear that both the Euler rule and DG-E are first order methods. Unexpectedly, DG-I turns out to be a second order method, even though the construction procedure has been identical. As said above, order analysis of discrete gradient methods is an interesting avenue for further research.

The picture changes radically when the step size is increased, even modestly to $h=7\cdot 10^{-4}$. To begin with, the trajectory computed by the Euler rule blow up to infinity, so it is not represented. Remarkably, the problem is not that of insufficient order: we tested an implicit Runge-Kutta method of order 2 (the basis of the \emph{ode23s} function in the \emph{Matlab ODE Suite}), and it also produced unbounded solutions. This is a significant finding, since methods designed for stiff differential equations are often assumed to reproduce better the qualitative behaviour, which is not the case here. Regarding the explicit discrete gradient method DG-E, its trajectory remains bounded, at least within the computed range, but the qualitative behaviour is completely wrong, as shown in Figure \ref{subfig:logistica_expl}. Instead of convergence to the equilibrium, undamped oscillations appear that destroy stability. In contrast, the correct behaviour is ultimately achieved by DG-I with the same step size, despite an initial transient, plotted in Figure \ref{subfig:logistica_impl}.

The apparent contradiction between the proved preservation of the Lyapunov function and the oscillatory solution provided by DG-E is explained by the local nature of the chosen Lyapunov function $V=\dfrac{1}{2}(1-y)^2$. The condition $\dfrac{dV}{dt}<0$ checked in Equation \eqref{eq:V-decreasing} only holds for $y<0$. This fact is dismissed in the original system, since the region $y<0$ cannot be reached from a positive initial value. However the discretization does take a step so large that the solution becomes negative. Once the function $V$ that drives the construction of the method is no longer a valid Lyapunov function, all bets are off. This suggests the first rule that must guide the construction of discrete gradient methods: \emph{find a Lyapunov function whose domain of validity is as large as possible}.

\section{Numerical experiments}
\label{sec:experiment}

In this section, we show the result of several numerical experiments designed to show the satisfactory performance of the designed discrete gradient method, assessed in terms of its ability to preserve the qualitative properties of the dynamical system. We are primarily interested in preserving the stability of the system, which will be evidenced by decreasing values of the considered Lyapunov function along solution trajectories of the numerical approximation. As a suitable case study, we propose the Duffing equation~\cite{calvo_projection_2010}, for which a Lyapunov function is known. The proposed method is compared with three conventional methods: the explicit Euler rule, a second-order Runge-Kutta method (RK2) that forms the basis of the \emph{ode23s} function in the \emph{Matlab ODE Suite}), and a fourth-order Runge-Kutta method (RK4), which the \emph{Matlab} \emph{ode45} function is based upon. Note that \emph{ode23s} is an implicit method, well suited to stiff equations, thus it is a strong competitor when preservation of qualitative features are considered, whereas \emph{ode45} is an explicit method design with higher order of accuracy in mind. In order to carry out a fair comparison among methods, all experiments are carried out with a fixed step size. Needless to say, our work on implementation of discrete gradient methods will eventually comprise variable step size mechanisms for error control.

All experiments apply to the Duffing equation that can be written as a first order system of ODEs $\frac{d\vec{y}}{dt}=f(\vec{y})$ by:
\begin{equation}
\label{eq:duffing}
\begin{aligned}
\frac{dy_{1}}{dt}&=y_{2} \\
\frac{dy_2}{dt}&=y_{1}-b\,y_{1}^3-a\,y_{2}
\end{aligned}
\end{equation}
with $b\neq 0$ and $a>0$. The system has three fixed points: $P_{0}=(0,0)$, $P_{1}=\left(\sqrt{1/b}, 0\right)$, and $P_{2}=\left(-\sqrt{1/b}, 0\right)$. A straightforward linearization shows  that $P_{0}$ is a saddle point, whereas $P_{1}$ and $P_{2}$ are stable equilibria. It is known that a Lyapunov function is defined by
\begin{equation}
    V(y_1,y_2)=\dfrac{1}{2}\,\left(y_{2}^2-y_{1}^2+\dfrac{b}{2}\,y_{1}^4\right)
\end{equation}
which has (local) minima at $P_1$ and $P_2$, since the gradient vanishes and the Hessian of $V$ is positive definite at both these points. The gradient of $V$ is the vector field:
\begin{equation}
\nabla V=
\begin{pmatrix}
-y_{1}+b\,y_1^3\\ 
y_2
\end{pmatrix}
\end{equation}
that leads to the energy-decreasing condition:
\begin{equation*}
\frac{dV}{dt}=\nabla V \cdot f(y)=-a\,y_{2}^2 \le 0   
\end{equation*}
Then, the system can be cast into the linear-gradient form, i.e.:
\begin{equation}
\frac{d\vec{y}}{dt}=
\begin{pmatrix}
0 & 1\\ 
-1 & -a
\end{pmatrix}
\begin{pmatrix}
-y_{1}+b\,y_1^3\\ 
y_2
\end{pmatrix}
\end{equation}
which entails the definition of the negative-definite matrix $L$:
\begin{equation}
L=
\begin{pmatrix}
0 & 1\\ 
-1 & -a
\end{pmatrix}
\end{equation}

Our implementation starts by computing the coordinate increment discrete gradient for the particular system given by Equation \eqref{eq:duffing}:
\begin{equation*}
    \overline{\nabla}V(y,z)=\frac{1}{2}\begin{pmatrix}
(z_{1}+y_{1})(-1+\dfrac{b}{2}(z_{1}^{2}+y_{1}^{2}))\\ 
z_{2}+y_{2}
\end{pmatrix}
\end{equation*} 
whereas we adopt the simplest approximation $\widetilde{L}=L$. Then, the discrete gradient method results:
\begin{equation*}
\begin{aligned}
z_{1} & =  y_{1}+\dfrac{h}{2}(z_{2}+y_{2})\\\\
z_{2} & =  \dfrac{h}{2}\left[-(z_{1}+y_{1})\left(-1+\dfrac{b}{2}(z_{1}^{2}+y_{1}^{2}) \right )-a(z_{2}+y_{2})\right ]
\end{aligned}
\end{equation*}
This implicit equation for $\vec{z}$ will be solved by Newton iteration until convergence at each time step.

\begin{figure} 
  \begin{subfigure}[b]{0.5\linewidth}
    \centering
    \includegraphics[width=0.75\linewidth]{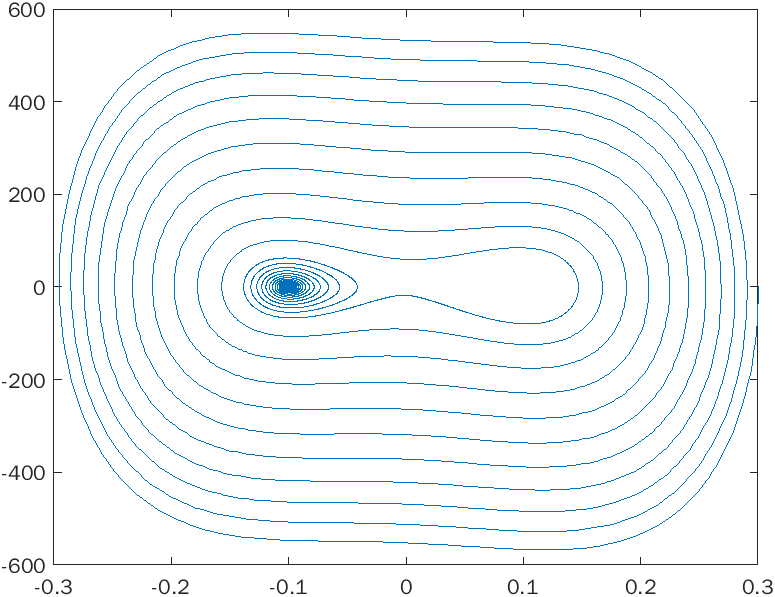} 
    \caption{Euler method} 
    \vspace{4ex}
  \end{subfigure}%% 
  \begin{subfigure}[b]{0.5\linewidth}
    \centering
    \includegraphics[width=0.75\linewidth]{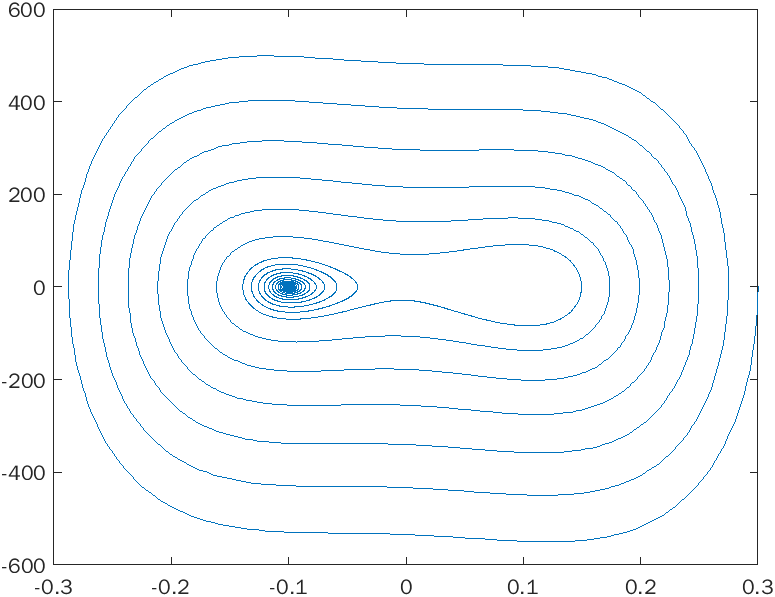} 
    \caption{Discrete gradient} 
    \vspace{4ex}
  \end{subfigure} 
  \begin{subfigure}[b]{0.5\linewidth}
    \centering
    \includegraphics[width=0.75\linewidth]{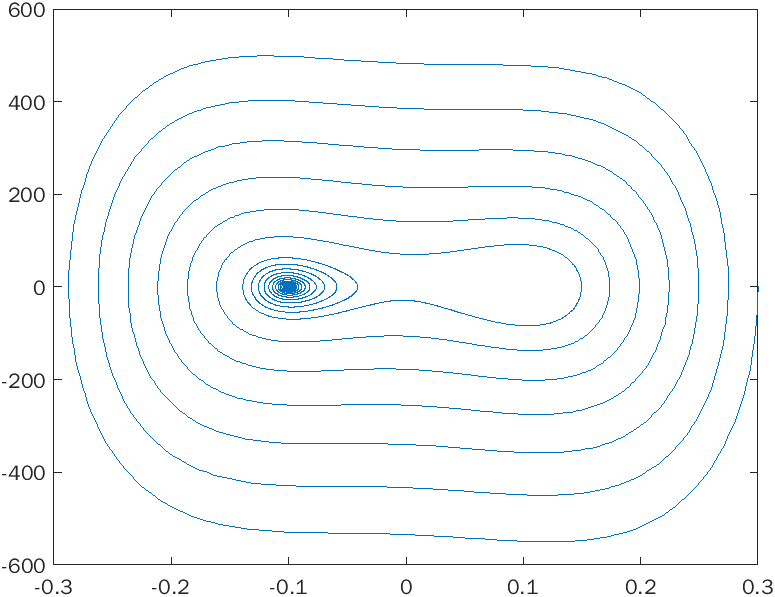} 
    \caption{RK of order 2} 
  \end{subfigure}%%
  \begin{subfigure}[b]{0.5\linewidth}
    \centering
    \includegraphics[width=0.75\linewidth]{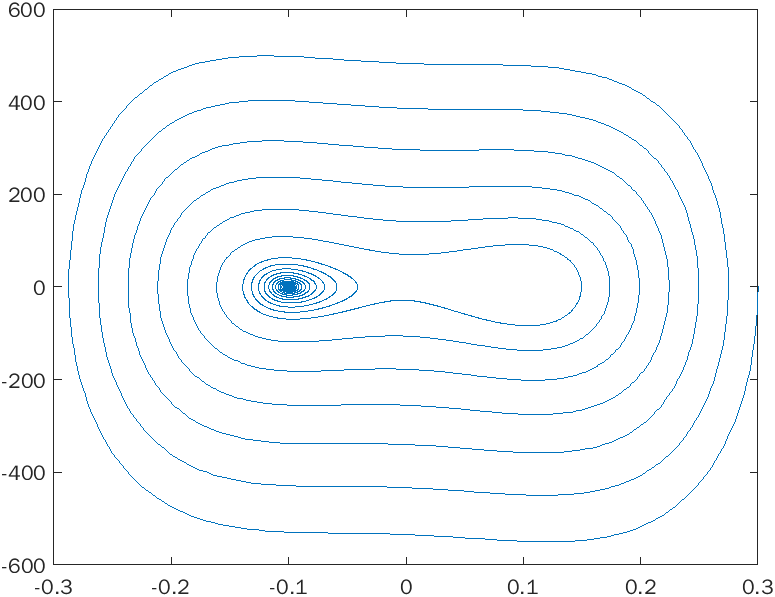} 
    \caption{RK of order 4}
    \end{subfigure} 
  \caption{Phase portrait for $h=10^{-5}$}
  \label{fig:1e-5} 
\end{figure}

\begin{figure} 
  \begin{subfigure}[b]{0.35\linewidth}
    \centering
    \includegraphics[width=0.75\linewidth]{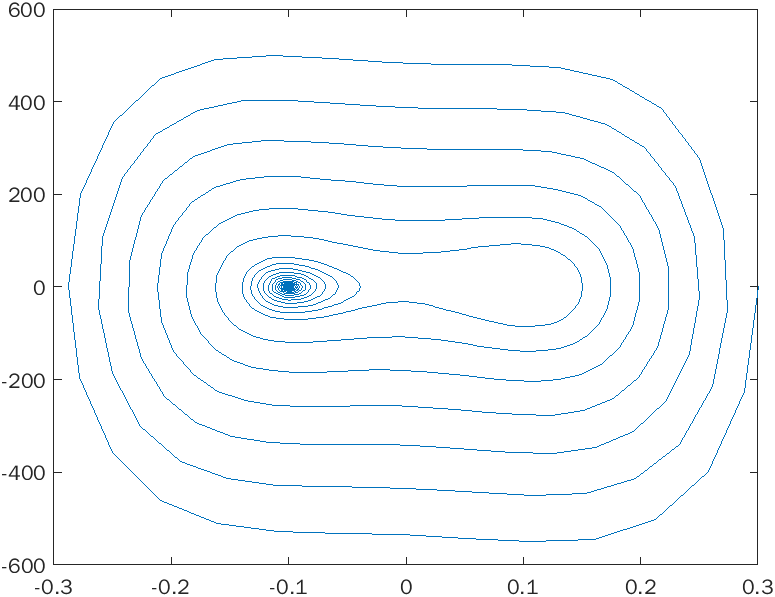} 
    \caption{Discrete gradient} 
  \end{subfigure} 
  \begin{subfigure}[b]{0.35\linewidth}
    \centering
    \includegraphics[width=0.75\linewidth]{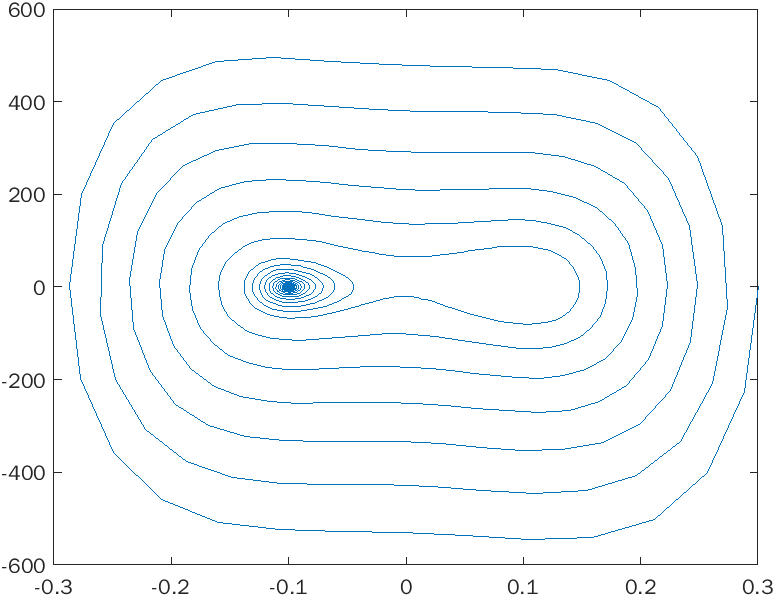} 
    \caption{RK of order 2} 
  \end{subfigure}%%
  \begin{subfigure}[b]{0.35\linewidth}
    \centering
    \includegraphics[width=0.75\linewidth]{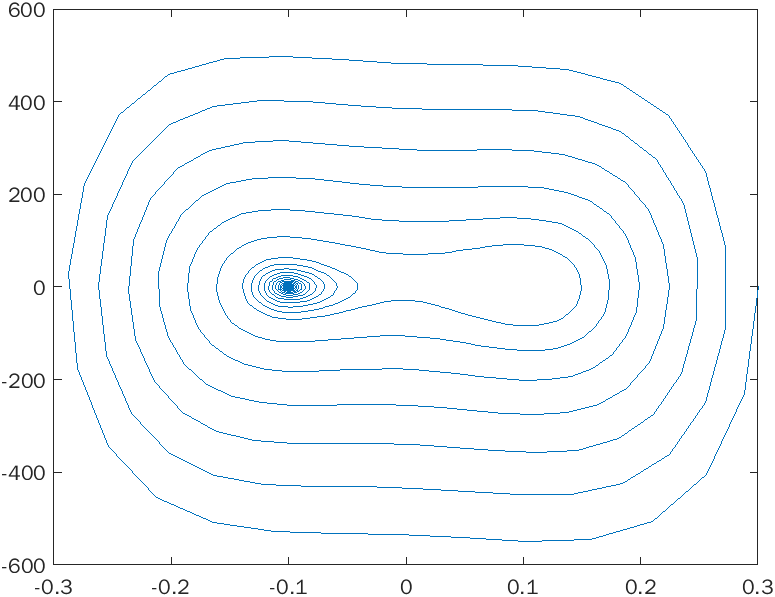} 
    \caption{RK of order 4}
\end{subfigure} 
  \caption{Phase portrait for $h=10^{-4}$}
  \label{fig:1e-4} 
\end{figure}

\begin{figure} 
  \begin{subfigure}[b]{0.5\linewidth}
    \centering
    \includegraphics[width=0.75\linewidth]{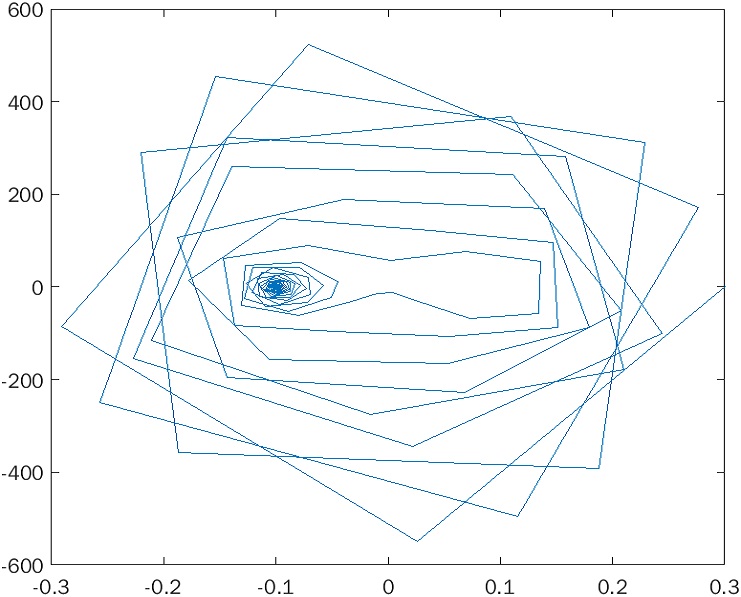} 
    \caption{Discrete gradient} 
    \label{GD-3} 
   % \vspace{4ex}
  \end{subfigure} 
  \begin{subfigure}[b]{0.5\linewidth}
    \centering
    \includegraphics[width=0.75\linewidth]{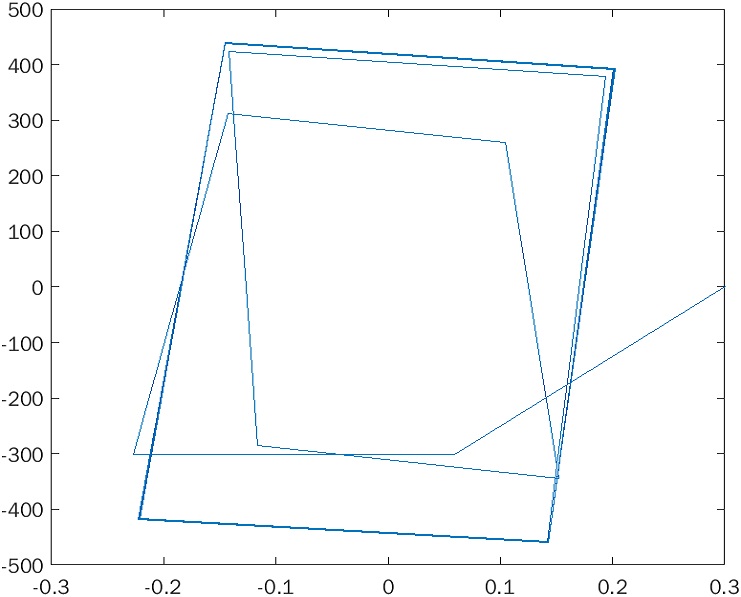} 
    \caption{RK of order 2} 
    \label{RK2-3} 
  \end{subfigure}
  \caption{Phase portrait for $h=10^{-3}$}
  \label{fig:1e-3} 
\end{figure}

All the experiments have been carried out considering $y_0=(0.3,0)$ as the initial point. We have designed three types of experiments. Firstly, we show the phase portrait that is obtained by applying each of the methods for different values of the  step size $h$ and compare it with the exact solution. Contrarily to the simple systems of the previous section, we do not have the benefit of an analytical solution, but we consider that the approximation obtained by the Euler's method with $h=10^{-8}$ is exact up to machine precision. The results of this set of experiments are shown in Figures \ref{fig:1e-5}, \ref{fig:1e-4} and \ref{fig:1e-3}. It can be seen how the behaviour of the discrete gradient method reproduces the phase portrait of the exact solution regardless of the step size. In contrast, the Euler's rule does not converge with step sizes greater than $10^{-5}$. As for the Runge-Kutta methods, both the order two and order four schemes fail when working with $h=10^{-3}$. Both explicit methods, Euler and RK4, produce trajectories that blow up towards unbounded values, thus they are not shown in figures. This is the case for both methods with $h=10^{-3}$ in Figure \ref{fig:1e-3} and the Euler's method with $h=10^{-4}$ in Figure \ref{fig:1e-4}. Despite the Euler's rule producing a bounded trajectory that converges to the stable equilibrium for small enough step size, the phase portrait is not correct. It is noticeable in Figure \ref{fig:1e-5}.a) that turns of the trajectory are closer than in other plots, suggesting that the numerical method is introducing a spurious dissipation.

\begin{table}
    \caption{Results of numerical experiments for the Duffing ODE.}
  \label{tab:rk}
	\centering
	\begin{tabular}{c|ccccc}
	\hline
Step size $h$ & Method & Comp. time & $\max \Delta V$
\\ \hline \hline
$ 10^{-3} $ & Euler & - &  $\infty$ \\\hline
        & RK4 & - & $ \infty $\\\hline
        & RK2 & 0.0146 & 0.0510  \\\hline
        & GD & 0.0069 &  $1.3010\cdot 10^{-18}$\\\hline \hline
        
$5\cdot 10^{-4}$ & Euler &  0.0027 & $\infty$ \\\hline 
        & RK4 & 0.0533 &  $8.6736\cdot 10^{-19}$\\\hline
        & RK2 &  0.0295 &  0.0091 \\\hline
        & GD & 0.0104 & $1.3010\cdot 10^{-18}$ \\\hline \hline
        
$ 10^{-4}$ & Euler & 0.0027 &  $\infty$\\\hline 
        & RK4 & 0.3558 & $8.6736\cdot 10^{-19}$ \\\hline
        & RK2 & 0.1299 &  $8.6736\cdot 10^{-19}$ \\\hline
        & GD &  0.0399 &  $1.3010\cdot 10^{-18}$\\\hline \hline
$5\cdot 10^{-5}$ & Euler & 0.0027 & $\infty$ \\\hline 
        & RK4 & 0.6885 &  $1.3010\cdot 10^{-18}$\\\hline
        & RK2 & 0.2429 &  $1.3010\cdot 10^{-18}$ \\\hline
        & GD & 0.1270 & $1.3010\cdot 10^{-18}$ \\\hline \hline  
$ 10^{-5}$ & Euler & 0.6185 & $2.8800\cdot 10^{-4}$ \\\hline 
        & RK4 & 3.2979  & $1.3010\cdot 10^{-18}$ \\\hline
        & RK2 & 1.2253 & $1.3010\cdot 10^{-18}$   \\\hline
        & GD & 0.5220 & $1.3010\cdot 10^{-18}$ \\\hline \hline 
$5\cdot 10^{-6}$ & Euler & 1.3717 & $7.2\cdot 10^{-5}$  \\\hline 
        & RK4 & 7.0309 & $1.3010\cdot 10^{-18}$ \\\hline
        & RK2 & 2.4288 &  $1.3010\cdot 10^{-18}$\\\hline
        & GD & 1.2380 & $1.3010\cdot 10^{-18}$ \\\hline \hline 
$ 10^{-6}$ & Euler & 3.4330 &  $2.88\cdot 10^{-6}$ \\\hline 
        & RK4 & 24.2176 & $1.7347\cdot 10^{-18}$ \\\hline
        & RK2 & 11.5226 &  $1.3010\cdot 10^{-18}$ \\\hline
        & GD & 4.7109 &  $1.7347\cdot 10^{-18}$ \\\hline           
\end{tabular}
\end{table}

\begin{figure} 
  \begin{subfigure}[b]{0.5\linewidth}
    \centering
    \includegraphics[width=0.75\linewidth]{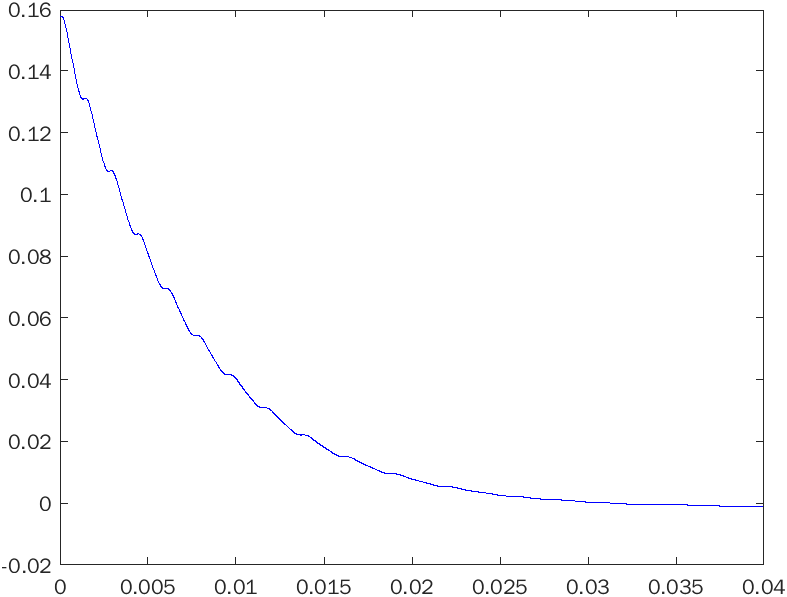} 
    \caption{Euler method} 
    \label{Euler-6} 
    \vspace{4ex}
  \end{subfigure}%% 
  \begin{subfigure}[b]{0.5\linewidth}
    \centering
    \includegraphics[width=0.75\linewidth]{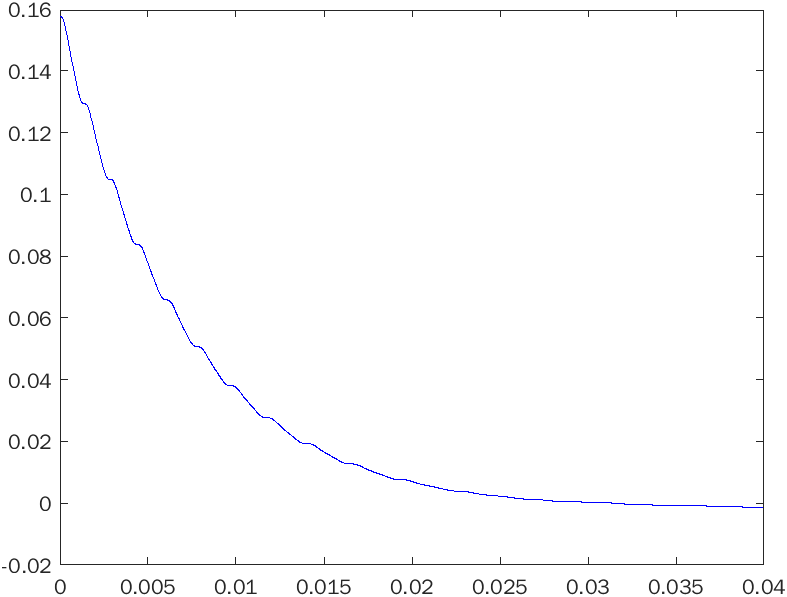} 
    \caption{Discrete gradient} 
    \label{GD-6} 
    \vspace{4ex}
  \end{subfigure} 
  \begin{subfigure}[b]{0.5\linewidth}
    \centering
    \includegraphics[width=0.75\linewidth]{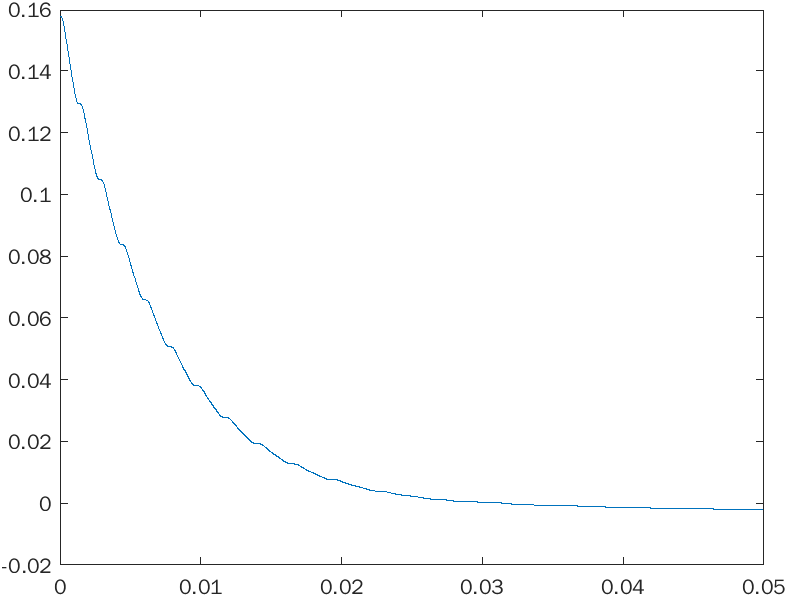} 
    \caption{RK of order two} 
    \label{RK2-6} 
  \end{subfigure}%%
  \begin{subfigure}[b]{0.5\linewidth}
    \centering
    \includegraphics[width=0.75\linewidth]{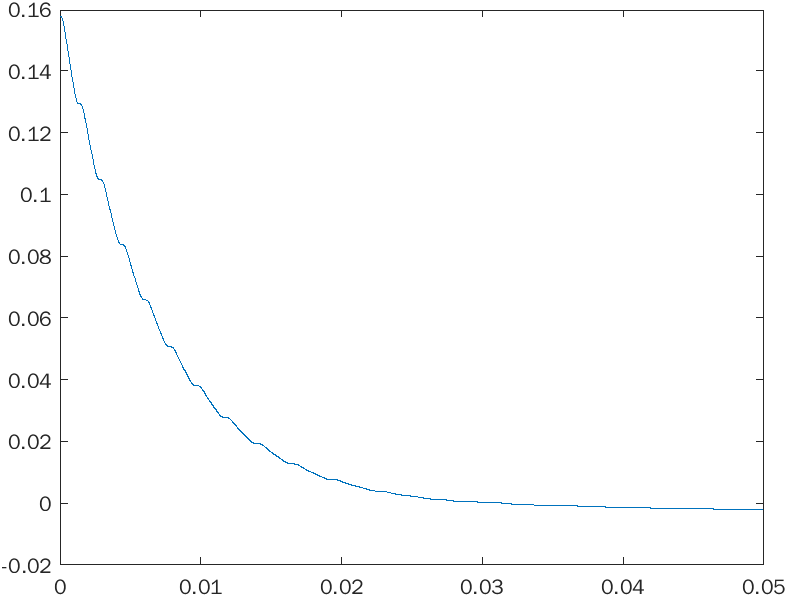} 
    \caption{RK of order four} 
    \label{RK4-6} 
  \end{subfigure} 
  \caption{Lyapunov function $h=10^{-6}$}
  \label{V-6} 
\end{figure}

\begin{figure} 
  \begin{subfigure}[b]{0.5\linewidth}
    \centering
    \includegraphics[width=0.75\linewidth]{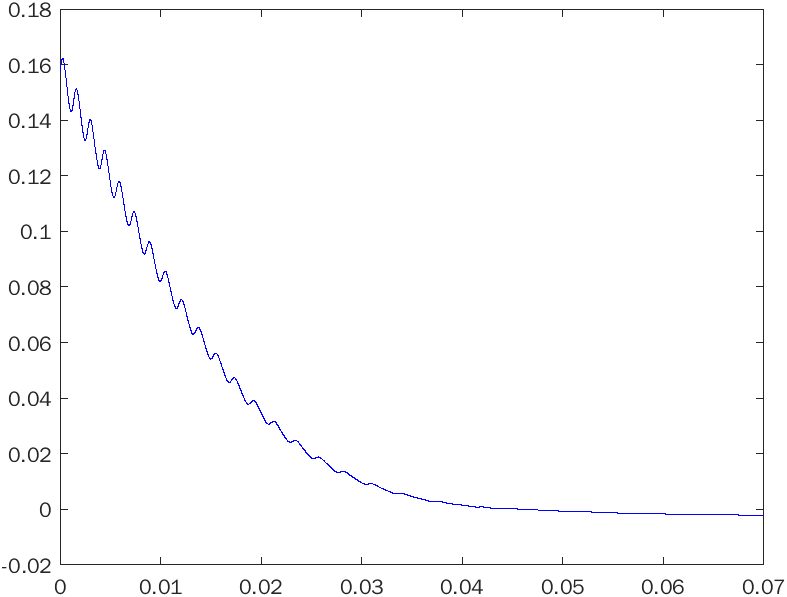} 
    \caption{Euler method} 
    \vspace{4ex}
  \end{subfigure}%% 
  \begin{subfigure}[b]{0.5\linewidth}
    \centering
    \includegraphics[width=0.75\linewidth]{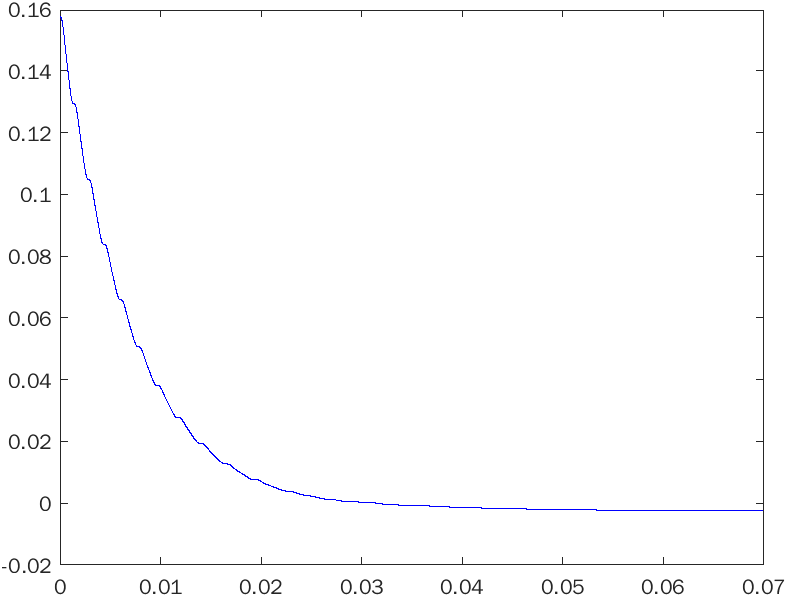} 
    \caption{Discrete gradient} 
    \vspace{4ex}
  \end{subfigure} 
  \begin{subfigure}[b]{0.5\linewidth}
    \centering
    \includegraphics[width=0.75\linewidth]{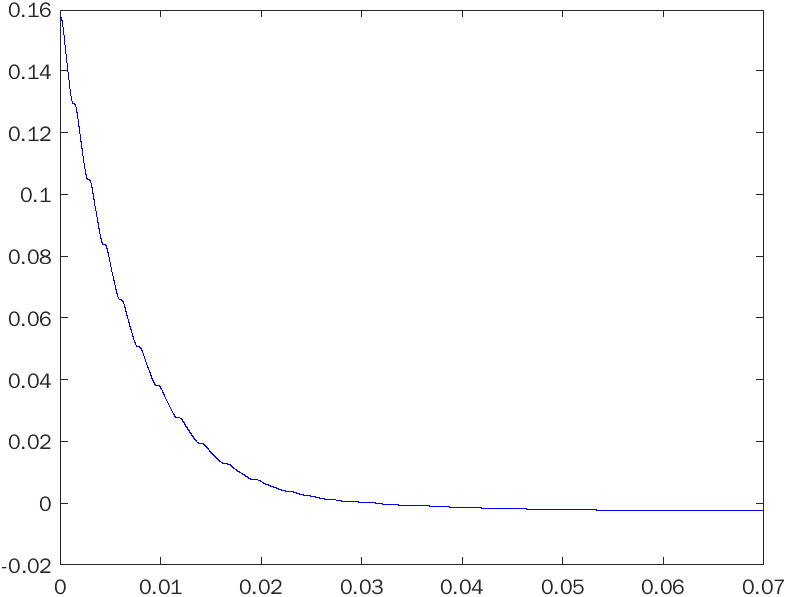} 
    \caption{RK of order two} 
  \end{subfigure}%%
  \begin{subfigure}[b]{0.5\linewidth}
    \centering
    \includegraphics[width=0.75\linewidth]{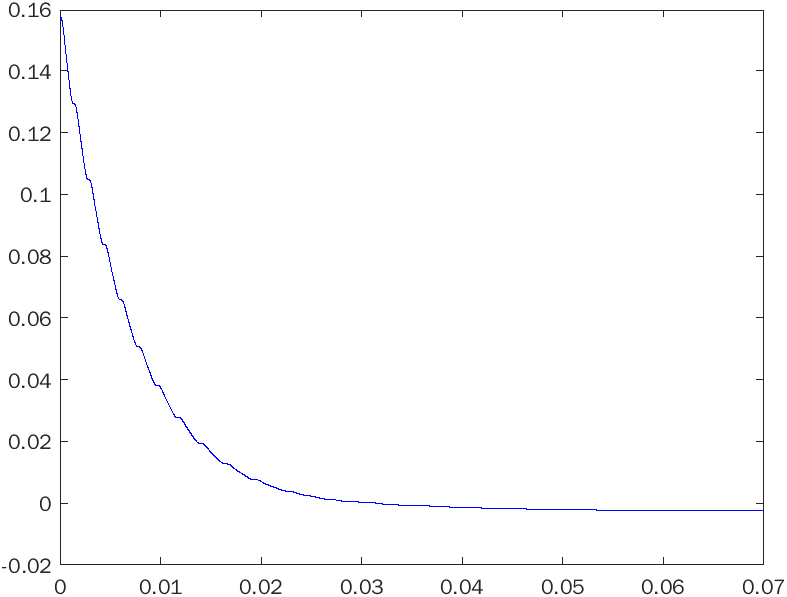} 
    \caption{RK of order four} 
  \end{subfigure} 
  \caption{Lyapunov function $h=10^{-5}$}
  \label{V-5} 
\end{figure}

\begin{figure} 
  \begin{subfigure}[b]{0.5\linewidth}
    \centering
    \includegraphics[width=0.75\linewidth]{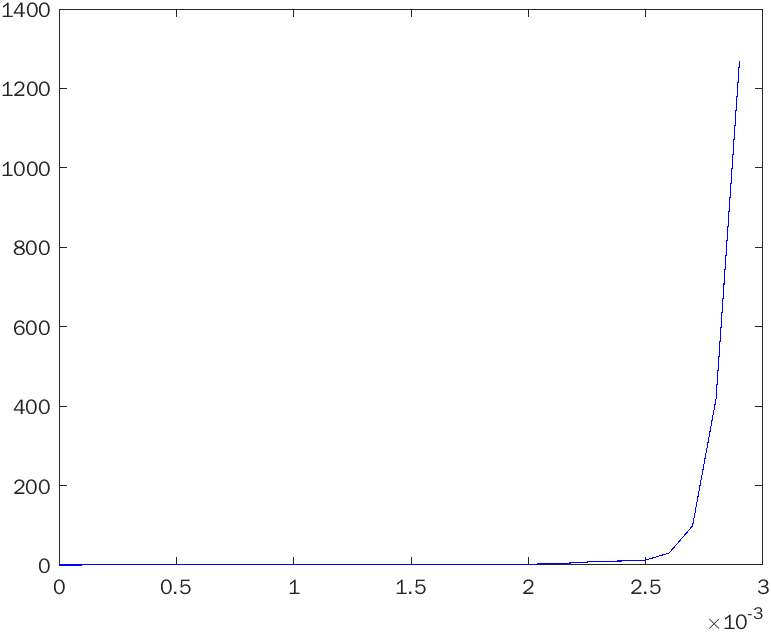} 
    \caption{Euler method} 
    \label{Euler-4} 
    \vspace{4ex}
  \end{subfigure}%% 
  \begin{subfigure}[b]{0.5\linewidth}
    \centering
    \includegraphics[width=0.75\linewidth]{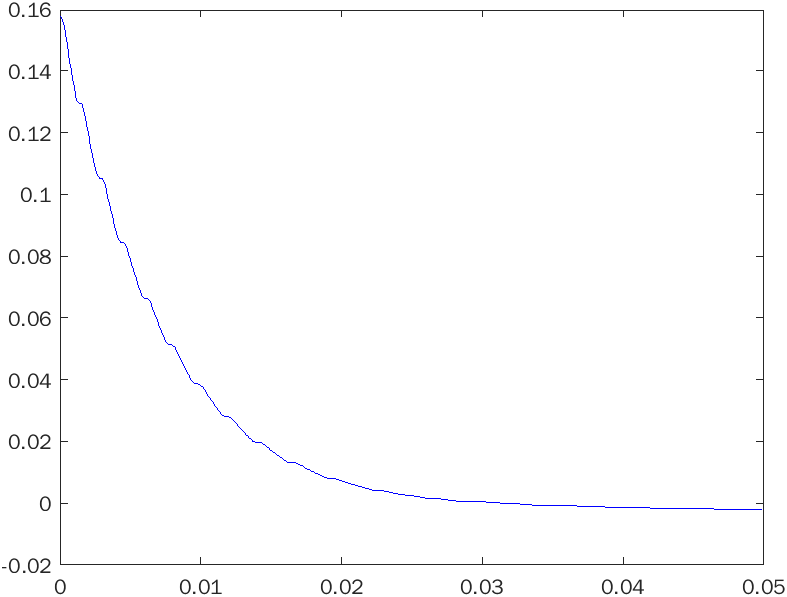} 
    \caption{Discrete gradient} 
    \vspace{4ex}
  \end{subfigure} 
  \begin{subfigure}[b]{0.5\linewidth}
    \centering
    \includegraphics[width=0.75\linewidth]{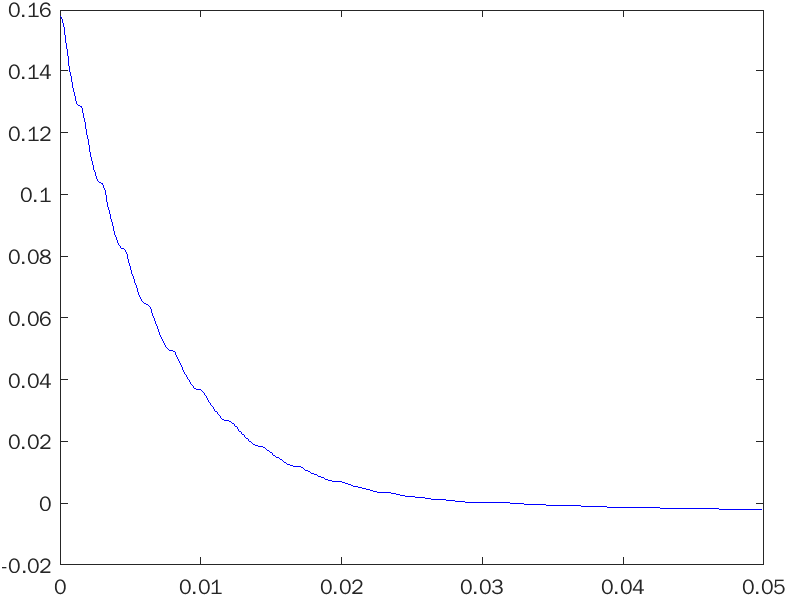} 
    \caption{RK of order two}  \end{subfigure}%%
  \begin{subfigure}[b]{0.5\linewidth}
    \centering
    \includegraphics[width=0.75\linewidth]{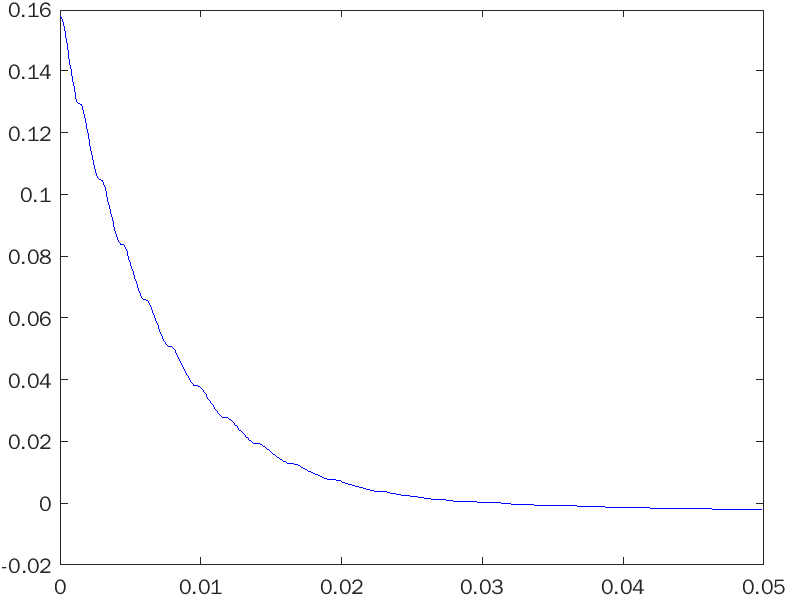} 
    \caption{RK of order four}  \end{subfigure} 
  \caption{Lyapunov function for $h=10^{-4}$}
  \label{V-4} 
\end{figure}

\begin{figure} 
  \begin{subfigure}[b]{0.5\linewidth}
    \centering
    \includegraphics[width=0.75\linewidth]{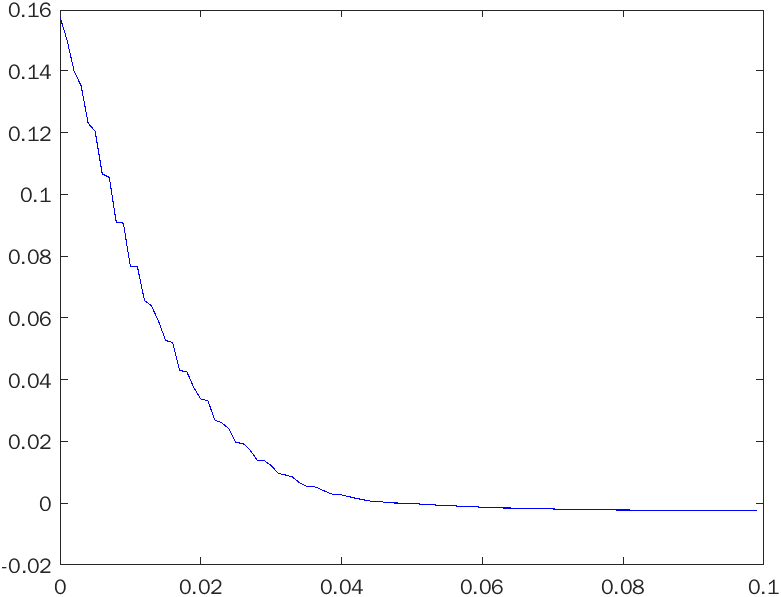} 
    \caption{Discrete gradient} 
    \label{V-GD-3} 
  \end{subfigure} 
  \begin{subfigure}[b]{0.5\linewidth}
    \centering
    \includegraphics[width=0.75\linewidth]{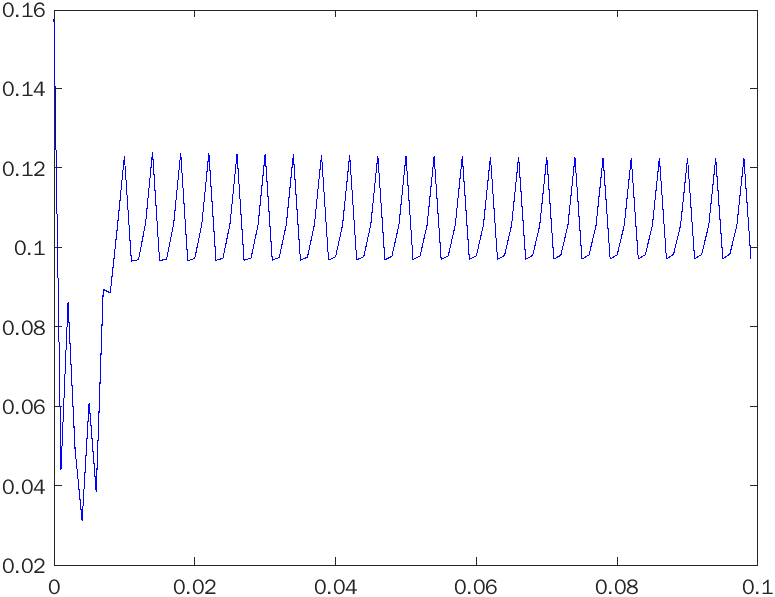} 
    \caption{RK of order two} 
    \label{V-RK2-3} 
  \end{subfigure}%%
  \caption{Lyapunov function for $h=10^{-3}$}
  \label{V-3} 
\end{figure}

On the other hand, taking into account that the fundamental objective of the designed method is the conservation of the Lyapunov function, we have designed another set of experiments focused on showing the behaviour of the Lyapunov function with respect to time. Table \ref{tab:rk} shows the values of the maximum increment of $V$ for each method and each step size used. We also plot in Figures \ref{V-6}-\ref{V-3} the trajectories of the value of $V$ for different step sizes. In general, it can be seen on the graphs that the Lyapunov function is decreasing along trajectories of the discrete gradient method, as expected by construction. The small positive increments shown in the table are within the range of machine precision so they are attributed to rounding rather than the numerical method. In contrast, much larger increases in the Lyapunov function are visible in Figure \ref{V-5} when using Euler's method with $h=10^{-5}$, even though for this step size the trajectories of the solution converge to the equilibrium. For large step sizes such as $h=10^{-3}$, only the implicit RK2 among conventional methods provides bounded trajectories. However, the evolution of $V$ shown in Figure \ref{V-3} reveals, even more clearly than the phase portrait, that the behaviour of the system is qualitatively corrupted. Periodic oscillations of $V$ prove that the system is not approaching an equilibrium and the Lyapunov function is no longer decreasing.

\begin{figure}
\centering
    \includegraphics[width=0.75\linewidth]{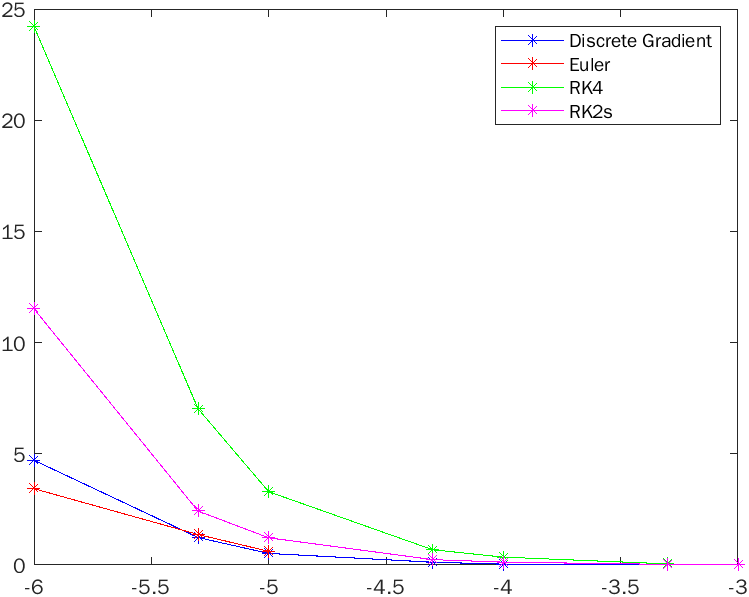} 
    \caption{Computational cost for different values of $h$} 
    \label{fig:coste} 
\end{figure}

Even when competitor conventional methods converge to a stable equilibrium, the proposed method is favourable in terms of computational cost. This is illustrated in Figure \ref{fig:coste} showing the real computation time for the different step sizes. The computing times are also shown  in Table \ref{tab:rk} for each combination of step size and method.

\section{Conclusions}
\label{sec:conclusion}
We have presented a methodology for implementation of numerical integrators that preserve a Lyapunov function of a dynamical system, namely discrete gradient methods. Analysis is performed on the proposed method, establishing that it is in principle a first order method, although the second order term is computed revealing the conditions for the method parameters under which a second order method would be obtained. As a proof of concept, a discrete gradient method is applied to the logistic equation, revealing the variety of choices that can lead to different numerical schemes with qualitatively different behaviours. The proposed method has been applied to the integration of the Duffing equation, which is regarded as a suitable test system: different parameter sets lead to oscillatory and stiff systems, whereas the preservation of the Lyapunov function is more important than accuracy of individual trajectories. Numerical experiments are also carried out to confirm the ability of discrete gradient methods to preserve the Lyapunov function, and the failure of standard Runge-Kutta codes for a wide range of step size values, since Lyapunov function increments occur, thus stability is lost.

We are currently engaged in further research in order to extend the results of this paper in several directions.  First, we are developing order conditions to obtain higher-order methods. Preliminary results show that this is possible, at least for order two, by defining the matrix $\widetilde{L}$ dependant not only on $\vec{y}$ and $\vec{z}$ but also on $h$. Another promising line considers composition and splitting techniques. The long-term objective would be to establish a systematic order theory for designing discrete gradient methods of arbitrary orders, in line with the recent paper \cite{Eidnes2022}. We are also trying to generalize the conditions for obtaining explicit methods, based on the original, implicit formulation.

This work suggests that general-purpose integrators are unable to keep pace with methods specifically designed for preserving the Lyapunov function. Thus we are extending our experiments in order to compare discrete gradient methods to both projection methods and Radau algorithms. In particular, it has been argued \cite{hairer_energy-diminishing_2013} that Radau methods are favourable due to its superior damping of high frequencies. In our experiments, we have detected that some discrete gradient methods possess an enhanced ability to deal with highly oscillatory systems. This question undoubtedly deserves deeper attention.  It also must be taken into account that the results of this paper are a proof of concept and much more can be done regarding implementation refinements of discrete gradient methods. The obvious advance is the inclusion of an error control device, which could derive from detecting lack of convergence of the Newton iteration. An improved discrete gradient method could be a serious competitor in applications where preserving the qualitative dynamical behaviour is more important that stringent accuracy of individual trajectories. For such systems, the integrators that preserve the Lyapunov function for arbitrary step sizes, such as discrete gradient methods, are endorsed as first line methods by our results.

\bibliographystyle{unsrtnat}
\bibliography{biblio.bib}

\end{document}